# Murnaghan-Nakayama Rules for Characters of Iwahori-Hecke Algebras of the Complex Reflection Groups $G(r,p,n)$

TOM HALVERSON
AND
ARUN RAM

Department of Mathematics
Macalester College
St. Paul, MN 55105

School of Mathematics and Statistics
University of Sydney
NSW 2006, Australia

ABSTRACT. Iwahori-Hecke algebras for the infinite series of complex reflection groups $G(r,p,n)$ were constructed recently in the work of Ariki and Koike, Broué and Malle, and Ariki. In this paper we give Murnaghan-Nakayama type formulas for computing the irreducible characters of these algebras. Our method is a generalization of that in our earlier paper in which we derived Murnaghan-Nakayama rules for the characters of the Iwahori-Hecke algebras of the classical Weyl groups. In both papers we have been motivated by C. Greene, who gave a new derivation of the Murnaghan-Nakayama formula for irreducible symmetric group characters by summing diagonal matrix entries in Young's seminormal representations. We use the analogous representations of the Iwahori-Hecke algebra of $G(r,p,n)$ given by Ariki.

## 1. INTRODUCTION

The finite irreducible complex reflection groups come in three infinite families: the symmetric groups $\mathcal{S}_n$ on $n$ letters; the wreath product groups $\mathbb{Z}_r \wr \mathcal{S}_n$, where $\mathbb{Z}_r$ denotes the cyclic group of order $r$; and a series of index-$p$ subgroups $G(r,p,n)$ of $\mathbb{Z}_r \wr \mathcal{S}_n$ for each positive integer $p$ that divides $r$. In the classification of finite irreducible reflection groups, besides these infinite families $\mathcal{S}_n$, $\mathbb{Z}_r$, and $G(r,p,n)$, there exist only 34 exceptional irreducible reflection groups, see [ST].

A formula for the irreducible characters of the Iwahori-Hecke algebras for $\mathcal{S}_n$ is known [Ram], [KW], [vdJ]. This formula is a $q$-analogue of the classical Murnaghan-Nakayama formula for computing the irreducible characters of $\mathcal{S}_n$. Similar formulas for the characters of the groups $G(r,p,n)$ are classically known, see [Mac], [Ste], [AK], [Osi] and the references there. Formulas of this type are also known for

1991 *Mathematics Subject Classification*. Primary 20C05; Secondary 05E05.
The second author gratefully acknowledges support from NSF grant DMS-9300523 and from a research fellowship under Australian Research Council grant No. A69330390.

Typeset by $\mathcal{A}_{\mathcal{M}}\mathcal{S}$-TEX

the Iwahori-Hecke algebras of Weyl groups of types B and D [HR], [Pfe1], [Pfe2]. Recently, Iwahori-Hecke algebras have been constructed for the groups $\mathbb{Z}_r \wr \mathcal{S}_n$ and $G(r, p, n)$ [AK], [BM], [Ari]. In this paper we derive Murnaghan-Nakayama type formulas for computing the irreducible characters of the Iwahori-Hecke algebras that correspond to $\mathbb{Z}_r \wr \mathcal{S}_n$ and $G(r, p, n)$.

Hoefsmit [Hfs] has given explicit analogues of Young's seminormal representations for the Iwahori-Hecke algebras of types $A_{n-1}$, $B_n$, and $D_n$. Ariki and Koike, [AK] and [Ari], have constructed "Hoefsmit-analogues" of Young's seminormal representations for Iwahori-Hecke algebras $H_{r,p,n}$ of the groups $G(r, p, n)$. Our approach is to derive the Murnaghan-Nakayama rules by computing the sum of diagonal matrix elements in an explicit "Hoefsmit" representation of each algebra. We are motivated by Curtis Greene [Gre], who takes this approach using the Young seminormal form of the irreducible representations of the symmetric group and gives a new derivation of the classical Murnaghan-Nakayama rule. Greene does this by using the Möbius function of a poset that is determined by the partition which indexes the irreducible representation. We generalize Greene's poset theorem so that it works for our cases. In this way we are able to compute the characters of the Hecke algebras $H_{r,n} = H_{r,1,n}$.

To compute the characters of the Iwahori-Hecke algebra $H_{r,p,n}$ of $G(r, p, n)$, $p > 1$, we use double centralizer methods (Clifford theory methods) to write these characters in terms of a certain bitrace on the irreducible representations of $H_{r,n} = H_{r,1,n}$. We then compute this bitrace in terms of the irreducible character values of $H_{r,n}$.

The character formulas given in this paper contain the Murnaghan-Nakayama rules for the complex reflection groups $G(r, p, n)$ and the Iwahori-Hecke algebras of classical type as special cases.

**Remark.** In this paper we only give formulas for computing the characters of certain "standard elements" of the Iwahori-Hecke algebra which are given by (2.10) in the case of the Iwahori-Hecke algebras of $\mathbb{Z}_r \wr \mathcal{S}_n$ and by (3.15) in the case of the Iwahori-Hecke algebras of $G(r, p, n)$, $p > 1$. In this paper we have not made any effort to show that this is sufficient to determine the values of the characters on all elements. We have a method for proving this which will be given in another paper. Results of this type for Iwahori-Hecke algebras of Weyl groups have been given in [GP].

## 2. Characters of Iwahori-Hecke Algebras of $(\mathbb{Z}/r\mathbb{Z}) \wr \mathcal{S}_n$

For positive integers $r$ and $n$, let $\mathcal{S}_n$ denote the symmetric group of order $n$ generated by $s_2, s_3, \ldots, s_n$, where $s_i$ denotes the transposition $s_i = (i-1, i)$, and let $\mathbb{Z}_r = \mathbb{Z}/r\mathbb{Z}$ denote the finite cyclic group of order $r$. Then the wreath product group $\mathbb{Z}_r \wr \mathcal{S}_n$ is a complex reflection group that can be identified with the group of all $n \times n$ permutation matrices whose non-zero entries are $r$th roots of unity.

Let $q$ and $u_1, u_2, \ldots, u_r$ be indeterminates. Let $H_{r,n}$ be the associative algebra with 1 over the field $\mathbb{C}(u_1, u_2, \ldots, u_r, q)$ given by generators $T_1, T_2, \ldots, T_n$ and

relations

(1) $T_i T_j = T_j T_i$,      for $|i - j| > 1$,
(2) $T_i T_{i+1} T_i = T_{i+1} T_i T_{i+1}$,      for $2 \leq i \leq n-1$,
(3) $T_1 T_2 T_1 T_2 = T_2 T_1 T_2 T_1$,
(4) $(T_1 - u_1)(T_1 - u_2) \cdots (T_1 - u_r) = 0$,
(5) $(T_i - q)(T_i + q^{-1}) = 0$,      for $2 \leq i \leq n$.

Upon setting $q = 1$ and $u_i = \xi^{i-1}$, where $\xi$ is a primitive $r$th root of unity, one obtains the group algebra $\mathbb{C}[\mathbb{Z}_r \wr \mathcal{S}_n]$. The algebras $H_{r,n}$ were first constructed by Ariki and Koike [AK], and they were classified as cyclotomic Hecke algebras of type $B_n$ by Broué and Malle [BM]. In the special case where $r = 1$ and $u_1 = 1$, we have $T_1 = 1$, and $H_{1,n}$ is isomorphic to an Iwahori-Hecke algebra of type $A_{n-1}$. The case $H_{2,n}$ when $r = 2$, $u_1 = p$, and $u_2 = p^{-1}$, is isomorphic to an Iwahori-Hecke algebra of type $B_n$.

**Shapes and Standard Tableaux.**

As in [Mac], we identify a partition $\alpha$ with its Ferrers diagram and say that a box $b$ in $\alpha$ is in position $(i,j)$ in $\alpha$ if $b$ is in row $i$ and column $j$ of $\alpha$. The rows and columns of $\alpha$ are labeled in the same way as for matrices.

An *$r$-partition* of size $n$ is an $r$-tuple, $\mu = (\mu^{(1)}, \mu^{(2)}, \ldots, \mu^{(r)})$ of partitions such that $|\mu^{(1)}| + |\mu^{(2)}| + \cdots + |\mu^{(r)}| = n$. If $\nu = (\nu^{(1)}, \nu^{(2)}, \ldots, \nu^{(r)})$ is another $r$-partition, we write $\nu \subseteq \mu$ if $\nu^{(i)} \subseteq \mu^{(i)}$ for $1 \leq i \leq r$. In this case, we say that $\mu/\nu = (\mu^{(1)}/\nu^{(1)}, \nu^{(2)}/\mu^{(2)}, \ldots, \mu^{(r)}/\nu^{(r)})$ is an $r$-skew shape. We refer to $r$-skew shapes and $r$-partitions collectively as *shapes*.

If $\lambda$ is a shape of size $n$, a *standard tableau* $L = (L^{(1)}, L^{(2)}, \ldots, L^{(r)})$ of shape $\lambda$ is a filling of the Ferrers diagram of $\lambda$ with the numbers $1, 2, \ldots, n$ such that the numbers are increasing left to right across the rows and increasing down the columns of each $L^{(i)}$. For any shape $\lambda$, let $\mathcal{L}(\lambda)$ denote the set of standard tableaux of shape $\lambda$ and, for each standard tableau $L$, let $L(k)$ denote the box containing $k$ in $L$.

**Representations.**

Define the *content* of a box $b$ of a (possibly skew) shape $\lambda = (\lambda^{(1)}, \ldots, \lambda^{(r)})$ by

$$\text{ct}(b) = u_k q^{2(j-i)}, \qquad \text{if } b \text{ is in position } (i,j) \text{ in } \lambda^{(k)}. \tag{2.1}$$

For each standard tableau $L$ of size $n$, define the scalar $(T_i)_{LL}$ by

$$(T_i)_{LL} = \frac{q - q^{-1}}{1 - \frac{\text{ct}(L(i-1))}{\text{ct}(L(i))}}, \qquad \text{for } 2 \leq i \leq n. \tag{2.2}$$

Note that $(T_i)_{LL}$ depends only on the *positions* of the boxes containing $i$ and $i-1$ in $L$.

Let $\lambda = (\lambda^{(1)}, \ldots, \lambda^{(r)})$ be a (possibly skew) shape of size $n$, and for each standard tableau $L \in \mathcal{L}(\lambda)$, let $v_L$ denote a vector indexed by $L$. Let $V^\lambda$ be the $\mathbb{C}(u_1, \ldots, u_r, q)$-vector space spanned by $\{v_L \mid L \in \mathcal{L}(\lambda)\}$, so that the vectors $v_L$ form a basis of $V^\lambda$. Define an action of $H_{r,n}$ on $V^\lambda$ by defining

$$\begin{aligned} T_1 v_L &= \text{ct}(L(1)) v_L, \\ T_i v_L &= (T_i)_{LL} v_L + (q^{-1} + (T_i)_{LL}) v_{s_i L}, \qquad 2 \leq i \leq n, \end{aligned} \tag{2.3}$$

where $s_i L$ is the same standard tableau as $L$ except that the positions of $i$ and $i-1$ are switched in $s_i L$. If $s_i L$ is not standard, then we define $v_{s_i L} = 0$.

The following theorem is due to Young [You] for the symmetric group $\mathcal{S}_n$, to Hoefsmit [Hfs] for $H_{1,n}$, and to Ariki and Koike [AK] for $H_{r,n}$, $r \geq 2$.

**Theorem 2.4.** ([You],[Hfs],[AK]) *The modules $V^\lambda$, where $\lambda$ runs over all $r$-partitions of size $n$, form a complete set of nonisomorphic irreducible modules for $H_{r,n}$.*

**Hoefsmit Elements.**

Define elements $t_i \in H_{r,n}$, for $1 \leq i \leq n$, by
$$t_i = T_i T_{i-1} \cdots T_2 T_1 T_2 \cdots T_{i-1} T_i. \tag{2.5}$$
To our knowledge, these elements were discovered by Hoefsmit in the case of the Iwahori-Hecke algebras of type $B_n$ and were rediscovered by Ariki and Koike for Iwahori-Hecke algebras of $\mathbb{Z}_r \wr \mathcal{S}_n$. For each standard tableau $L$ of size $n$, define the scalar $(t_i)_{LL}$ by
$$(t_i)_{LL} = \mathrm{ct}(L(i)), \qquad \text{for } 1 \leq i \leq n. \tag{2.6}$$
The following proposition is due to Hoefsmit for $r = 1, 2$ and to Ariki and Koike for $r > 2$.

**Proposition 2.7.** ([Hfs], Prop. 3.3.3; [AK], Prop. 3.16) *For $1 \leq i \leq n$, the action of $t_i$ on a vector $v_L$, where $L$ is a standard tableau, is*
$$t_i v_L = (t_i)_{LL} v_L.$$
Furthermore, these elements commute:

**Proposition 2.8.** ([AK], Lemma 3.3) *The subalgebra $\mathfrak{U}_{n,r}$ of $H_{r,n}$ generated by $t_1, t_2, \ldots, t_n$ is an abelian subalgebra, i.e., the $t_i$ commute.*

**Standard Elements.**

For $1 \leq k < \ell \leq n$ and $0 \leq i \leq r-1$, define
$$R_{k\ell}^{(i)} = (t_k)^i T_{k+1} T_{k+2} \cdots T_\ell \tag{2.9}$$
and, for each $1 \leq k \leq n$, define $R_{kk}^{(i)} = (t_k)^i$. We say that an $\mathcal{S}_n$-*sequence* of length $m$ is a sequence $\vec{\ell} = (\ell_1, \ldots, \ell_m)$ satisfying $1 \leq \ell_1 < \ell_2 < \cdots < \ell_m = n$, and we say that a $\mathbb{Z}_r$-*sequence* of length $m$ is a sequence $\vec{i} = (i_1, \ldots, i_m)$ satisfying $0 \leq i_j \leq r-1$ for each $j$. For an $\mathcal{S}_n$-sequence $\vec{\ell} = (\ell_1, \ldots, \ell_m)$ and a $\mathbb{Z}_r$-sequence $\vec{i} = (i_1, \ldots, i_m)$, define
$$T_{\vec{\ell}}^{\vec{i}} = R_{1,\ell_1}^{(i_1)} R_{\ell_1+1,\ell_2}^{(i_2)} \cdots R_{\ell_{m-1}+1,\ell_m}^{(i_m)} \in H_{r,n}. \tag{2.10}$$
For example, in $H_{4,10}$ we have the standard element
$$T_{(3,4,8,10)}^{(0,2,3,1)} = R_{1,3}^{(0)} R_{4,4}^{(2)} R_{5,8}^{(3)} R_{9,10}^{(1)} = T_2 T_3 (t_4)^2 (t_5)^3 T_6 T_7 T_8 t_9 T_{10}.$$

For $1 \leq k < \ell \leq n$ and $0 \leq i \leq r-1$, we define
$$\Delta_{k\ell}^{(i)}(L) = (t_k)_{LL}^i (T_{k+1})_{LL} (T_{k+2})_{LL} \cdots (T_\ell)_{LL}, \tag{2.11}$$
and for $1 \leq k \leq n$, we define $\Delta_{kk}^{(i)}(L) = (t_k)_{LL}^i$. Since $(T_j)_{LL}$ depends only on the positions of the boxes $j$ and $j-1$ in $L$, the scalar $\Delta_{k\ell}^{(i)}(L)$ depends only on the positions of the boxes containing $k, k+1, \ldots, \ell$ in $L$.

**Proposition 2.12.** *Let $\vec{\ell} = (\ell_1, \ldots, \ell_m)$ be an $\mathcal{S}_n$-sequence and $\vec{i} = (i_1, \ldots, i_m)$ be a $\mathbb{Z}_r$-sequence, and let $L$ be a standard tableau of size $n$. Let $T_{\vec{\ell}}^{\vec{i}} v_L \big|_{v_L}$ denote the coefficient of $v_L$ in $T_{\vec{\ell}}^{\vec{i}} v_L$. Then*

$$T_{\vec{\ell}}^{\vec{i}} v_L \big|_{v_L} = \Delta_{1,\ell_1}^{(i_1)}(L) \Delta_{\ell_1+1,\ell_2}^{(i_2)}(L) \cdots \Delta_{\ell_{m-1}+1,\ell_m}^{(i_m)}(L).$$

*In particular, for given sequences $\vec{\ell}$ and $\vec{i}$, the value $T_{\vec{\ell}}^{\vec{i}} v_L \big|_{v_L}$ depends only on the positions and the linear order of the boxes in $L$.*

*Proof.* This follows from the definition of the action of $H_{r,n}$ on standard tableaux and the fact (2.3) that when $T_i$ acts on a standard tableau $L$ it affects only the positions of $L$ containing $i$ and $i-1$. The result follows, since $t_j$ acts as a scalar (Prop. 2.7), and $T_{\vec{\ell}}^{\vec{i}}$ otherwise is a product (from right to left) of a decreasing sequence of generators $T_i$. □

**Characters.**

If $L$ is a standard tableau (of any shape, possibly of skew shape) with $n$ boxes, define

$$\Delta^{(i)}(L) = \Delta_{1,n}^{(i)}(L), \tag{2.13}$$

and for any shape $\lambda$ (possibly skew), define

$$\Delta^{(i)}(\lambda) = \sum_{L \in \mathcal{L}(\lambda)} \Delta^{(i)}(L). \tag{2.14}$$

In making these definitions, the actual values in the boxes of $L$ do not matter, only their positions and their order relative to one another are relevant. Thus, the definitions make sense when the standard tableaux have values that form a subset of $\{1, 2, \ldots\}$ (with the usual linear order).

For an $r$-partition $\lambda$, let $\chi_{H_{r,n}}^\lambda$ denote the character of the irreducible $H_{r,n}$-representation $V^\lambda$ determined in Theorem 2.4. The following theorem is our analogue of the Murnaghan-Nakayama rule.

**Theorem 2.15.** *Let $\vec{\ell} = (\ell_1, \ldots, \ell_m)$ be an $\mathcal{S}_n$-sequence, $\vec{i} = (i_1, \ldots, i_m)$ be a $\mathbb{Z}_r$-sequence, and suppose that $\lambda$ is an $r$-partition of size $n$. Then*

$$\chi_{H_{r,n}}^\lambda(T_{\vec{\ell}}^{\vec{i}}) = \sum_{\emptyset = \mu^{(0)} \subseteq \mu^{(1)} \subseteq \cdots \subseteq \mu^{(m)} = \lambda} \Delta^{(i_1)}(\mu^{(1)}) \Delta^{(i_2)}(\mu^{(2)}/\mu^{(1)}) \cdots \Delta^{(i_m)}(\mu^{(m)}/\mu^{(m-1)}),$$

*where the sum is over all sequences of shapes $\emptyset = \mu^{(0)} \subseteq \mu^{(1)} \subseteq \cdots \subseteq \mu^{(m)} = \lambda$ such that $|\mu^{(j)}/\mu^{(j-1)}| = |\ell_j|$.*

*Proof.* By Proposition 2.12 the character $\chi_{H_{r,n}}^\lambda$ is given by

$$\chi_{H_{r,n}}^\lambda(T_{\vec{\ell}}^{\vec{i}}) = \sum_{L \in \mathcal{L}(\lambda)} T_{\vec{\ell}}^{\vec{i}} v_L \big|_{v_L} = \sum_{L \in \mathcal{L}(\lambda)} \Delta_{1,\ell_1}^{(i_1)}(L) \Delta_{\ell_1+1,\ell_2}^{(i_2)}(L) \cdots \Delta_{\ell_{m-1}+1,\ell_m}^{(i_m)}(L).$$

The result follows by collecting terms according to the positions occupied by the various segments of the numbers $\{1, 2, \ldots, \ell_1\}$, $\{\ell_1 + 1, \ldots, \ell_2\}$, $\ldots$, $\{\ell_{m-1} + 1, \ldots, \ell_m\}$. □

In view of Theorem 2.15 it is desirable to give an explicit formula for the value of $\Delta^{(i)}(\lambda)$. To do so requires some further notations: The shape $\lambda$ is a *border strip* if it is connected and does not contain two boxes which are adjacent in the same northwest-to-southeast diagonal. This is equivalent to saying that $\lambda$ is connected and does not contain any $2 \times 2$ block of boxes. The shape $\lambda$ is a *broken border strip* if it does not contain any $2 \times 2$ block of boxes. Therefore, a broken border strip is a union of connected components, each of which is a border strip.

Drawing Ferrers diagrams as in [Mac], we say that a *sharp corner* in a border strip is a box with no box above it and no box to its left. A *dull corner* in a border strip is a box that has a box to its left and a box above it but has no box directly northwest of it. The picture below shows a broken border strip with two connected components where each of the sharp corners has been marked with an **s** and each of the dull corners has been marked with a **d**.

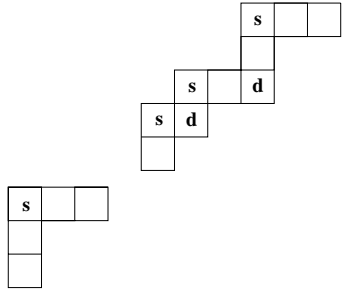

Figure 2.16

The following theorem is proved using Corollary 4.14 of Theorem 4.6. We have placed these results in Section 4, because they stand on their own as results on planar posets.

**Theorem 2.17.** *Let $\lambda$ be any shape (possibly skew) with $n$ boxes. Let $\mathrm{CC}$ be the set of connected components of $\lambda$, and let $cc = |\mathrm{CC}|$ be the number of connected components of $\lambda$.*

*(a) If $\lambda$ is not a broken border strip, then $\Delta^{(k)}(\lambda) = 0$;*

*(b) If $\lambda$ is a broken border strip, then*

$$\Delta^{(0)}(\lambda) = (q - q^{-1})^{cc-1} \prod_{bs \in \mathrm{CC}} q^{c(bs)-1}(-q^{-1})^{r(bs)-1},$$

*and, for $1 \leq k \leq r - 1$,*

$$\Delta^{(k)}(\lambda) = (-q + q^{-1})^{cc-1} \left( \prod_{s \in \mathrm{SC}} \mathrm{ct}(s) \right) \left( \prod_{d \in \mathrm{DC}} \mathrm{ct}(d)^{-1} \right)$$
$$\times \sum_{t=0}^{|\mathrm{DC}|} (-1)^t e_t(\mathrm{ct}(\mathrm{DC})) h_{k-t-cc}(\mathrm{ct}(\mathrm{SC}))$$
$$\times \prod_{bs \in \mathrm{CC}} q^{c(bs)-1}(-q^{-1})^{r(bs)-1},$$

*where $\mathrm{SC}$ and $\mathrm{DC}$ denote the set of sharp corners and dull corners in $\lambda$, respectively, and if $bs$ is a border strip, then $r(bs)$ is the number of rows in $bs$, and $c(bs)$ is*

the number of columns in $bs$. The content $\operatorname{ct}(b)$ of a box $b$ is as given in (2.1). The function $e_t(\operatorname{ct}(DC))$ is the elementary symmetric function in the variables $\{\operatorname{ct}(d), d \in DC\}$, and the function $h_{k-t-cc}(\operatorname{ct}(SC))$ is the homogeneous symmetric function in the variables $\{\operatorname{ct}(s), s \in SC\}$.

*Proof.* Recall from (2.2) that

$$(T_k)_{LL} = \frac{q - q^{-1}}{1 - \frac{\operatorname{ct}(L(k-1))}{\operatorname{ct}(L(k))}}.$$

It follows from the definitions of $\Delta^{(k)}(L)$ in (2.13) and (2.14) that we may apply Corollary 4.14 with $x_b = \operatorname{ct}(b)$ for all boxes $b$ in $\lambda$.

For two boxes $a$ and $b$ in $\lambda$ that are adjacent in a diagonal, we have

$$\frac{1 - \operatorname{ct}(a)\operatorname{ct}(b)^{-1}}{(q - q^{-1})} = \frac{1 - 1}{q - q^{-1}} = 0.$$

Thus, $\Delta^{(k)}(\lambda) = 0$ if $\lambda$ is not a broken border strip. Furthermore,

$$\frac{(q - q^{-1})}{1 - \operatorname{ct}(a)\operatorname{ct}(b)^{-1}} = \begin{cases} \frac{(q - q^{-1})}{1 - q^{-2}} = q, & \text{if } a \text{ and } b \text{ are adjacent in a row,} \\ \frac{(q - q^{-1})}{1 - q^2} = -q^{-1}, & \text{if } a \text{ and } b \text{ are adjacent in a column.} \end{cases}$$

The result now follows from Corollary 4.14. $\square$

## 3. Characters of Iwahori-Hecke Algebras of $G(r,p,n)$

In this section we define the complex reflection groups $G(r,p,n)$ and their Iwahori-Hecke algebras $H_{r,p,n}$. The groups $G(r,p,n)$ are normal subgroups of index $p$ in the groups $G(r,1,n)$, and the groups $G(r,1,n)$ are isomorphic to the wreath products $\mathbb{Z}_r \wr S_n$. The corresponding Hecke algebras $H_{r,p,n}$ are subalgebras of $H_{r,n}$. We compute the irreducible characters of $H_{r,p,n}$ in terms of the irreducible characters of $H_{r,n}$, which are computed in Section 2.

**The Complex Reflection Groups $G(r,p,n)$.**

Let $r, p, d$, and $n$ be positive integers such that $pd = r$. The complex reflection group $G(r,p,n)$ is the set of $n \times n$ matrices such that

(a) The entries are either 0 or $r$th roots of unity.
(b) There is exactly one nonzero entry in each row and each column.
(c) The $d$th power of the product of the nonzero entries is 1.

The order of $G(r,p,n)$ is given by $|G(r,p,n)| = dr^{n-1}n!$, and $G(r,p,n)$ is a normal subgroup of $G(r,1,n)$ of index $p$.

Let $\zeta = e^{2\pi i/r}$ be a primitive $r$th root of unity. Then $G(r,p,n)$ is generated by the elements

$$s_0 = \zeta^p E_{11} + \sum_{i=2}^n E_{ii}, \qquad s_1 = \zeta E_{12} + \zeta^{-1} E_{21} + \sum_{i=3}^n E_{ii},$$

$$s_j = \sum_{i \neq j, j-1} E_{ii} + E_{(j-1)j} + E_{j(j-1)}, \qquad 2 \leq j \leq n,$$

where $E_{ij}$ denotes the $n \times n$ matrix with a 1 in the $i$th row and $j$th column and with all other entries 0.

**Example 3.1.** The following are important special cases of $G(r,p,n)$.
 (1) $G(1,1,n) = S_n$, the symmetric group.
 (2) $G(r,1,n) = \mathbb{Z}_r \wr S_n$.
 (3) $G(2,1,n) = WB_n$ the Weyl group of type B.
 (4) $G(2,2,n) = WD_n$ the Weyl group of type D.

**The Hecke algebras.**

Let $\varepsilon = e^{2\pi i/p}$ be a primitive $p$th root of unity, and let $q$ and $x_0^{1/p}, \ldots, x_{d-1}^{1/p}$ be indeterminates. Then $H_{r,n}$ is the associative algebra with 1 over the field $\mathbb{C}(x_0^{1/p}, \ldots, x_{d-1}^{1/p}, q)$ given by generators $T_1, \ldots, T_n$, and relations

 (1) $T_i T_j = T_j T_i$,      for $|i-j| > 1$,
 (2) $T_i T_{i+1} T_i = T_{i+1} T_i T_{i+1}$,      for $2 \le i \le n-1$,
 (3) $T_1 T_2 T_1 T_2 = T_2 T_1 T_2 T_1$,
 (4) $(T_1^p - x_0)(T_1^p - x_1) \cdots (T_1^p - x_{d-1}) = 0$,
 (5) $(T_i - q)(T_i + q^{-1}) = 0$,      for $2 \le i \le n$.

This is the same as the definition of the algebra $H_{r,n}$ in section 2 except that we are using $\varepsilon^\ell x_k^{1/p}$, $0 \le k \le d-1$, $0 \le \ell \le p-1$, in place of $u_1, \ldots, u_r$. Let $H_{r,p,n}$ be the subalgebra of $H_{r,n}$ generated by the elements

$$a_0 = T_1^p, \quad a_1 = T_1^{-1} T_2 T_1, \quad \text{and} \quad a_i = T_i, \quad 2 \le i \le n. \tag{3.2}$$

Ariki ([Ari], Proposition 1.6) shows that $H_{r,p,n}$ is an analogue of the Iwahori-Hecke algebra for the groups $G(r,p,n)$. The special case $H_{2,2,n}$ is isomorphic to an Iwahori-Hecke algebra of type $D_n$.

**Shapes and Tableaux.**

As above, $r, p, d$, and $n$ are positive integers such that $pd = r$. We organize each $r$-partition $\lambda$ of size $n$ into $d$ groups of $p$ partitions each, so that we can write

$$\lambda = (\lambda^{(k,\ell)}), \quad \text{for } 0 \le k \le d-1 \text{ and } 0 \le \ell \le p-1,$$

where each $\lambda^{(k,\ell)}$ is a partition and $\sum_{k,\ell} |\lambda^{(k,\ell)}| = n$. It is convenient to view the partitions $\lambda^{(k,0)}, \ldots, \lambda^{(k,p-1)}$ as all lying on a circle so that we have $d$ necklaces of partitions, each necklace with $p$ partitions on it. In order to specify this arrangement, we shall say that $\lambda$ is a $(d,p)$-partition.

As in section 2, we let $\mathcal{L}(\lambda)$ denote the set of standard tableaux of shape $\lambda$, and, for each standard tableau $L$, let $L(i)$ denote the box containing $i$ in $L$.

**Action on Standard Tableaux.**

Let $\lambda = (\lambda^{(k,\ell)})$ be a $(d,p)$-partition of size $n$. Since $H_{r,p,n}$ is a subalgebra of $H_{r,n}$, the irreducible $H_{r,n}$-representations $V^\lambda$ are (not necessarily irreducible) representations of $H_{r,p,n}$. However we can easily describe the action of $H_{r,p,n}$ on $V^\lambda$ by restricting the action of $H_{r,n}$.

With the given specializations of the $u_i$, the content of a box $b$ of $\lambda$, see (2.1), is

$$\mathrm{ct}(b) = \varepsilon^\ell x_k^{1/p} q^{2(j-i)}, \qquad \text{if } b \text{ is in position } (i,j) \text{ in } \lambda^{(k,\ell)}.$$

As in Section 2, we define, for each standard tableau $L$ of size $n$, the scalar

$$(T_i)_{LL} = \frac{q - q^{-1}}{1 - \frac{\mathrm{ct}(L(i-1))}{\mathrm{ct}(L(i))}}, \qquad \text{for } 2 \leq i \leq n.$$

From (2.3) and (3.2), it follows that the action of $H_{r,p,n}$ on $V^\lambda$ is given by

$$\begin{aligned}
a_0 v_L &= \mathrm{ct}(L(1))^p v_L = x_k v_L, \quad \text{if } 1 \in L^{(k,\ell)}, \\
a_1 v_L &= (T_2)_{LL} v_L + \frac{\mathrm{ct}(L(1))}{\mathrm{ct}(s_2 L(1))}(q^{-1} + (T_2)_{LL}) v_{s_2 L}, \\
a_i v_L &= (T_i)_{LL} v_L + (q^{-1} + (T_i)_{LL}) v_{s_i L}.
\end{aligned} \qquad (3.3)$$

Recall that $t_i = T_i \cdots T_2 T_1 T_2 \cdots T_i$ for $1 \leq i \leq n$, and define elements $S_i \in H_{r,p,n}$, $1 \leq i \leq n$, by

$$\begin{aligned}
S_1 &= a_0 = t_1^p, \\
S_2 &= a_1 a_2 = t_1^{-1} t_2, \\
S_i &= a_i a_{i-1} \cdots a_4 a_3 a_1 a_2 a_3 a_4 \cdots a_{i-1} a_i = t_1^{-1} t_i, \qquad \text{for } 3 \leq i \leq n.
\end{aligned}$$

It follows from the action of the $t_i$ (Prop. 2.7) that the action of $S_i$ on $V^\lambda$ is also diagonal and is given by

$$\begin{aligned}
S_1 v_L &= \mathrm{ct}(L(1))^p v_L \\
S_i v_L &= \mathrm{ct}(L(1))^{-1} \mathrm{ct}(L(i)) v_L.
\end{aligned} \qquad (3.4)$$

Furthermore, since the $t_i$ commute (Prop. 2.8), it follows that the $S_i$ commute.

**A $\mathbb{Z}/p\mathbb{Z}$ action on shapes.**

Let $\lambda = (\lambda^{(k,\ell)})$ be a $(d,p)$-partition. We define an operation $\sigma$ that moves the partitions on each circle over one position. Given a box $b$ in position $(i,j)$ of the partition $\lambda^{(k,\ell)}$ then $\sigma(b)$ is the same box $b$ except moved to be in position $(i,j)$ of $\lambda^{(k,\ell+1)}$, where $\ell+1$ is taken modulo $p$. The map $\sigma$ is an operation of order $p$ and acts uniformly on the shape $\lambda = (\lambda^{(k,\ell)})$, on standard tableaux $L = (L^{(k,\ell)})$ of shape $\lambda$, and on the basis vector $v_L$ of $V^\lambda$:

$$\sigma(\lambda) = (\lambda^{(k,\ell+1)}), \quad \sigma(L) = (L^{(k,\ell+1)}), \quad \text{and} \quad \sigma(v_L) = v_{\sigma(L)}.$$

We use the notation $\sigma$ in each case, since the operation is always clear from context. In the last case, extend linearly to get the vector space homomorphism $\sigma: V^\lambda \longrightarrow V^{\sigma(\lambda)}$. If $b$ is a box in a shape $\lambda$, then

$$\mathrm{ct}(\sigma(b)) = \varepsilon \mathrm{ct}(b). \qquad (3.5)$$

**Lemma 3.6.** *The map $\sigma: V^\lambda \to V^{\sigma(\lambda)}$ is a $H_{r,p,n}$-module isomorphism, i.e., $\sigma$ commutes with the action of $H_{r,p,n}$.*

*Proof.* Since $(T_i)_{LL}$, see (2.2), depends only on the row and column of boxes $i$ and $i-1$ and not on the position of the tableaux, we have $(T_i)_{\sigma L, \sigma L} = (T_i)_{LL}$, for all $1 \leq i \leq n$ and all standard tableaux $L$. Since $\sigma(s_2 L) = s_2 \sigma(L)$, it follows that $a_i v_{\sigma L} = \sigma(a_i v_L)$. (Note that, because $\varepsilon$ is a $p$th root of unity, $T_1$ does *not* commute with $\sigma$ but that $a_0$ does.) □

The set of transformations

$$\{\sigma^\alpha \mid 0 \leq \alpha \leq p-1\}$$

defines an action of the cyclic group $\mathbb{Z}/p\mathbb{Z}$ on the set of $(d,p)$-partitions and on the set of vector spaces $V^\lambda$.

**Irreducible Representations.**

Fix a $(d,p)$-partition $\lambda$ of size $n$, and let $K_\lambda$ be the stabilizer of $\lambda$ under the action of $\mathbb{Z}/p\mathbb{Z}$. The group $K_\lambda$ is a cyclic group of order $|K_\lambda|$ and is generated by the transformation $\sigma^{f_\lambda}$ where $f_\lambda$ is the smallest integer between 1 and $p$ such that $\sigma^{f_\lambda}(\lambda) = \lambda$. Thus,

$$K_\lambda = \{\sigma^{\alpha f_\lambda} : V^\lambda \to V^\lambda \mid 0 \leq \alpha \leq |K_\lambda| - 1\}. \tag{3.7}$$

Figure (3.9) is an example of a (3,6)-partition $\lambda$ for which $f_\lambda = 2$ and $K_\lambda = \{1, \sigma^2, \sigma^4\} \cong \mathbb{Z}_3$. The elements of $K_\lambda$ are all $H_{r,p,n}$-module isomorphisms. The irreducible $K_\lambda$-modules are all one-dimensional, and the characters of these modules are given explicitly by

$$\eta_j : \begin{array}{ccc} K_\lambda & \longrightarrow & \mathbb{C} \\ \sigma^{f_\lambda} & \longmapsto & \varepsilon^{j f_\lambda} \end{array} \tag{3.8}$$

where $0 \leq j \leq |K_\lambda| - 1$. To see this note that $\omega = \varepsilon^{f_\lambda}$ is a primitive $|K_\lambda|$-th root of unity.

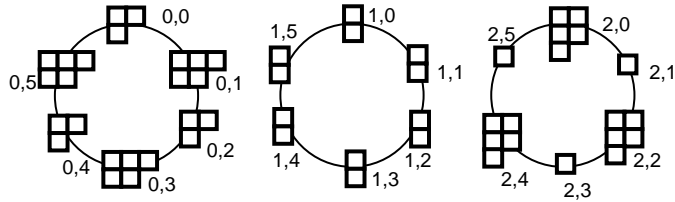

Figure 3.9. A (3,6)-partition $\lambda$ with $f_\lambda = 2$

It follows (from a standard double centralizer result) that as an $H_{r,p,n} \times K_\lambda$-bimodule

$$V^\lambda \cong \bigoplus_{j=0}^{|K_\lambda|-1} V^{(\lambda,j)} \otimes Z_j, \tag{3.10}$$

where $V^{(\lambda,j)}$ is an $H_{r,p,n}$-module and $Z_j$ is the irreducible $K_\lambda$-module with character $\eta_j$. Ariki ([Ari], Theorem 2.6) has explicitly constructed the modules $V^{(\lambda,j)}$ and proved that they form a complete set of irreducible $H_{r,p,n}$-modules. From the point of view of (3.10), one can prove that the $V^{(\lambda,j)}$ are irreducible $H_{r,p,n}$-modules by setting $q = 1$ and $x_k = 1$ for all $0 \leq k \leq d-1$ and appealing to the corresponding result for the group $G(r,p,n)$.

**Theorem 3.11.** ([Ari], Theorem 2.6) *The modules $V^{(\lambda,j)}$, where $\lambda$ runs over all $r$-partitions and $0 \le j \le |K_\lambda| - 1$, form a complete set of nonisomorphic irreducible modules for $H_{r,p,n}$.*

**Remark 3.12.** It should be noted that if $f_\lambda = p$ and thus $|K_\lambda| = 1$, then the irreducible $H_{r,n}$-module $V^\lambda$ is also an irreducible $H_{r,p,n}$-module.

### Characters.

Fix a $(d,p)$-partition $\lambda$ and let $\chi^{(\lambda,j)}$ denote the character of the irreducible $H_{r,p,n}$-module $V^{(\lambda,j)}$ defined by (3.10). Let $\chi^\lambda$ denote the $H_{r,p,n} \times K_\lambda$-bitrace on the module $V^\lambda$, i.e. if $h \in H_{r,p,n}$ and $\sigma^{\alpha f_\lambda} \in K_\lambda$, then

$$\chi^\lambda(h\sigma^{\alpha f_\lambda}) = \sum_{L \in \mathcal{L}(\lambda)} h\sigma^{\alpha f_\lambda} v_L \Big|_{v_L} = \sum_{L \in \mathcal{L}(\lambda)} hv_L \Big|_{v_{\sigma^{-\alpha f_\lambda} L}}, \qquad (3.13)$$

where $h\sigma^{\alpha f_\lambda} v_L \big|_{v_L}$ denotes the coefficient of $v_L$ in the expansion of $h\sigma^{\alpha f_\lambda} v_L$ in terms of the basis of $V^\lambda$ corresponding to standard tableaux.

By taking traces in the module equation (3.10), we obtain

$$\chi^\lambda(h\sigma^{\alpha f_\lambda}) = \sum_{j=0}^{|K_\lambda|-1} \chi^{(\lambda,j)}(h) \eta_j(\sigma^{\alpha f_\lambda}) = \sum_{j=0}^{|K_\lambda|-1} \chi^{(\lambda,j)}(h) \varepsilon^{j \alpha f_\lambda}.$$

By the orthogonality of characters for $K_\lambda$ (or by direct computation) this formula can be inverted to give

$$\chi^{(\lambda,j)}(h) = \frac{1}{|K_\lambda|} \sum_{\alpha=0}^{|K_\lambda|-1} \varepsilon^{-j\alpha f_\lambda} \chi^\lambda(h\sigma^{\alpha f_\lambda}), \quad \text{where } f_\lambda = p/|K_\lambda|. \qquad (3.14)$$

### Standard elements.

For $1 \le k \le n$, define $S^{(i)}_{kk} = S^i_k$ and define $\tilde{S}^{(i)}_{12} = S^i_1 a_1$. For all other $k < \ell$, define

$$S^{(i)}_{k\ell} = S^i_k a_{k+1} \cdots a_\ell, \quad \text{and} \quad \tilde{S}^{(i)}_{1\ell} = S^i_1 a_1 a_3 \cdots a_\ell. \qquad (3.15)$$

Following the definitions in (2.10), let $(\ell_1, \ldots, \ell_m)$ be an $\mathcal{S}_n$-sequence and let $(i_1, \ldots, i_m)$ be a $\mathbb{Z}_r$-sequence. The remainder of this section is devoted to computing the values

$$\chi^{(\lambda,j)}(\tilde{S}^{i_1}_{1\ell_1} S^{(i_2)}_{\ell_1+1,\ell_2} \cdots S^{(i_m)}_{\ell_{m-1}+1,\ell_m}) \quad \text{and} \quad \chi^{(\lambda,j)}(S^{i_1}_{1\ell_1} S^{(i_2)}_{\ell_1+1,\ell_2} \cdots S^{(i_m)}_{\ell_{m-1}+1,\ell_m}).$$

### Reduction to $R^{(i_1)}_{1,\ell_1} R^{(i_2)}_{\ell_1+1,\ell_2} \cdots R^{(i_m)}_{\ell_{m-1}+1,\ell_m}$.

Recall the definition (2.9) of the element $R^{(i)}_{m,\ell}$ of $H_{r,n}$. We now show that it is sufficient to compute characters on special products of these elements.

**Lemma 3.16.** *The group generated by $\{a_i \mid 0 \le i \le n\}$ in $H_{r,n}$ is a normal subgroup of the group generated by $\{T_j \mid 1 \le j \le n\}$ in $H_{r,n}$.*

*Proof.* It is sufficient to show that $T_j a_i T_j^{-1}$ and $T_j^{-1} a_i T_j$ can be written as a product of the $a_i$s and their inverses. The only two nontrivial calculations are the following.

$$T_1 a_2 T_1^{-1} = T_1 T_2 T_1^{-1} = T_2^{-1} T_2 T_1 T_2 T_1^{-1} = T_2^{-1} T_1^{-1} T_2 T_1 T_2 = a_2^{-1} a_1 a_2$$

and
$$T_1^{-1}a_1T_1 = T_1^{-2}T_2T_1^2 = T_1^{-1}(T_1^{-1}T_2T_1T_2)T_2^{-1}T_1$$
$$= T_1^{-1}T_2T_1T_2T_1^{-1}T_2^{-1}T_1 = a_1a_2a_1^{-1}. \qquad \square$$

**Lemma 3.17.** *For a word $g = T_{i_1}^{\pm 1}\cdots T_{i_m}^{\pm 1}$ in the generators of $H_{r,n}$ (i.e., an element of the group in $H_{r,n}$ generated by the $T_i$), let $\beta(g)$ denote the number of $T_1$'s minus the number of $T_1^{-1}$'s in $g$ so that $\beta(g)$ is the net number of $T_1$'s in the word. Let $h$ be in the group generated by $\{a_i \mid 0 \le i \le n\}$ in $H_{r,p,n}$. Then*
$$\chi^\lambda(ghg^{-1}\sigma^{\alpha f_\lambda}) = \varepsilon^{-\alpha f_\lambda \beta(g)}\chi^\lambda(h\sigma^{\alpha f_\lambda}).$$

*Proof.* First note that, by Lemma 3.16, $ghg^{-1} \in H_{r,p,n}$, so it makes sense to consider the bitrace. Then
$$\chi^\lambda(ghg^{-1}\sigma^{\alpha f_\lambda}) = \chi^\lambda(hg^{-1}\sigma^{\alpha f_\lambda}g).$$
We must be very careful here, because, although the action of $h$ commutes with the action of $\sigma^{\alpha f_\lambda}$, the action of $g$ and $g^{-1}$ do not. In fact, since
(1) $T_1 v_{\sigma L} = \mathrm{ct}(\sigma L(1))v_{\sigma L} = \varepsilon \mathrm{ct}(L(1))v_{\sigma L} = \varepsilon \sigma(T_1 v_L)$, and
(2) $T_i v_{\sigma L} = \sigma(T_i v_L)$, for $2 \le i \le n$, by Lemma 3.6, and
it follows that $g^{-1}\sigma^{\alpha f_\lambda} = \varepsilon^{-\alpha f_\lambda \beta(g)}\sigma^{\alpha f_\lambda}g^{-1}$. Thus,
$$\chi^\lambda(ghg^{-1}\sigma^{\alpha f_\lambda}) = \varepsilon^{-\alpha f_\lambda \beta(g)}\chi^\lambda(gh\sigma^{\alpha f_\lambda}g^{-1}) = \varepsilon^{-\alpha f_\lambda \beta(g)}\chi^\lambda(h\sigma^{\alpha f_\lambda}). \qquad \square$$

**Lemma 3.18.** *Let $(\ell_1,\ldots,\ell_m)$ be an $\mathcal{S}_n$-sequence, and let $(i_1,\ldots,i_m)$ be a $\mathbb{Z}_r$-sequence. satisfying $0 \le i_j \le r-1$ for each $j$.*
$$\chi^\lambda(\tilde{S}_{1\ell_1}^{i_1}S_{\ell_1+1,\ell_2}^{(i_2)}\cdots S_{\ell_{m-1}+1,\ell_m}^{(i_m)}\sigma^{\alpha f_\lambda})$$
$$= \varepsilon^{\alpha f_\lambda(i_2+\cdots+i_m-1)}\chi^\lambda(R_{1,\ell_1}^{(i_1 p - i_2 - \cdots - i_m)}R_{\ell_1+1,\ell_2}^{(i_2)}\cdots R_{\ell_{m-1}+1,\ell_m}^{(i_m)}\sigma^{\alpha f_\lambda}),$$
$$\chi^\lambda(S_{1\ell_1}^{i_1}S_{\ell_1+1,\ell_2}^{(i_2)}\cdots S_{\ell_{m-1}+1,\ell_m}^{(i_m)}\sigma^{\alpha f_\lambda})$$
$$= \varepsilon^{\alpha f_\lambda(i_2+\cdots+i_m)}\chi^\lambda(R_{1,\ell_1}^{(i_1 p - i_2 - \cdots - i_m)}R_{\ell_1+1,\ell_2}^{(i_2)}\cdots R_{\ell_{m-1}+1,\ell_m}^{(i_m)}\sigma^{\alpha f_\lambda}).$$

*Proof.* We have
$$S_{k\ell}^{(i)} = t_1^{-i}R_{k\ell}^{(i)}, \quad \text{and} \quad S_{1\ell}^{(i)} = t_1^{ip-1}T_2 t_1 T_3 \cdots T_\ell.$$
Since $t_1$ commutes with $T_i$ for $i > 2$ it follows that, for any $\mathcal{S}_n$-sequence $(\ell_1,\ldots,\ell_m)$ and $\mathbb{Z}_r$-sequence $(i_1,\ldots,i_m)$, we have
$$\tilde{S}_{1\ell_1}^{i_1}S_{\ell_1+1,\ell_2}^{(i_2)}\cdots S_{\ell_{m-1}+1,\ell_m}^{(i_m)} = t_1^{i_1 p - 1}T_2 t_1^{1-i_2-\cdots-i_m}T_3\cdots T_{\ell_1}R_{\ell_1+1,\ell_2}^{(i_2)}\cdots R_{\ell_{m-1}+1,\ell_m}^{(i_m)}$$
$$S_{1\ell_1}^{i_1}S_{\ell_1+1,\ell_2}^{(i_2)}\cdots S_{\ell_{m-1}+1,\ell_m}^{(i_m)} = t_1^{i_1 p}T_2 t_1^{-i_2-\cdots-i_m}T_3\cdots T_{\ell_1}R_{\ell_1+1,\ell_2}^{(i_2)}\cdots R_{\ell_{m-1}+1,\ell_m}^{(i_m)}.$$
Both of these can be conjugated by a power of $t_1$ to give
$$R_{1,\ell_1}^{(i_1 p - i_2 - \cdots - i_m)}R_{\ell_1+1,\ell_2}^{(i_2)}\cdots R_{\ell_{m-1}+1,\ell_m}^{(i_m)}.$$
Now use Lemma 3.17. $\square$

In view of Lemma 3.18 and (3.14) we shall try to compute the values of $\chi^\lambda(h\sigma^{\alpha f_\lambda})$, for elements $h \in H_{r,p,n}$ of the form $R_{1,\ell_1}^{(i_1)}R_{\ell_1+1,\ell_2}^{(i_2)}\cdots R_{\ell_{m-1}+1,\ell_m}^{(i_m)}$, where $i_1 + \cdots + i_m \equiv 0 \pmod{p}$, and where $\vec{\ell} = (\ell_1,\ldots,\ell_m)$ is an $\mathcal{S}_n$-sequence and $\vec{i} = (i_1,\ldots,i_m)$ is a $\mathbb{Z}_r$-sequence. In fact we shall prove the following theorem. We state the theorem now in order to establish the notations.

**Theorem 3.19.** Let $\lambda$ be a $(d,p)$-partition, where $pd = r$. Let $\alpha$ be such that $0 \leq \alpha \leq |K_\lambda| - 1$ where $K_\lambda$ is as defined in (3.7). Define

$$f_\lambda = p/|K_\lambda|, \quad \text{and} \quad \gamma = \frac{|K_\lambda|}{\gcd(\alpha, |K_\lambda|)}.$$

Let

$$h = R^{(i_1)}_{1,\ell_1} R^{(i_2)}_{\ell_1+1,\ell_2} \cdots R^{(i_m)}_{\ell_{m-1}+1,\ell_m}$$

where $(\ell_1, \cdots, \ell_m)$ is an $\mathcal{S}_n$-sequence and $(i_1, \cdots, i_m)$ is a $\mathbb{Z}_r$-sequence such that $i_1 + \cdots + i_m = 0 \pmod{p}$. The element $h$ is an element of $H_{r,p,n} \subseteq H_{r,n}$. If all $\ell_i$ in the sequence $(\ell_1, \ldots, \ell_m)$ are divisible by $\gamma$ then define

$$\bar{n} = n/\gamma, \quad \bar{r} = r/\gamma, \quad \bar{p} = p/\gamma,$$
$$(\bar{\ell}_1, \ldots, \bar{\ell}_m) = (\ell_1/\gamma, \ldots, \ell_m/\gamma),$$
$$\bar{\lambda}^{(k,\tau)} = \lambda^{(k,\tau)}, \quad \text{for } 0 \leq \tau \leq \bar{p} - 1 \quad \text{and} \quad \bar{h} = R^{(0)}_{1,\bar{\ell}_1} \cdots R^{(0)}_{\bar{\ell}_{m-1}+1,\bar{\ell}_m}.$$

Then:
(a) If $\ell_i$ is not divisible by $\gamma$ for some $1 \leq i \leq m$ then $\chi^\lambda(h\sigma^{\alpha f_\lambda}) = 0$.

(b) If all $\ell_i$ are divisible by $\gamma$ and if $i_k \neq 0$ for some $k$, then $\chi^\lambda(h\sigma^{\alpha f_\lambda}) = 0$.

(c) If all $\ell_i$ are divisible by $\gamma$ and if $i_k = 0$ for all $k$, then

$$\chi^\lambda(h\sigma^{\alpha f_\lambda}) = \frac{\gamma^{\bar{n}}}{[\gamma]^{\bar{n}-m}} \chi^{\bar{\lambda}}_{H_{\bar{r},\bar{n}}}(\bar{h}) \prod_{i=1}^{\gamma} \left( \frac{q}{1-\varepsilon^{-i}} + \frac{q^{-1}}{1-\varepsilon^i} \right)^{\bar{n}},$$

where $H_{\bar{r},\bar{n}}$ is with parameter $q^\gamma$, in place of $q$ and with parameters $\varepsilon^{\gamma\tau} x_k^\gamma$, $0 \leq k \leq d-1$, $0 \leq \tau \leq \bar{p} - 1\}$ in place of $u_1, \ldots, u_{\bar{r}}$. The element $\bar{h}$ is viewed as an element of the algebra $H_{\bar{r},\bar{n}}$ and $[\gamma] = (q^\gamma - q^{-\gamma})/(q - q^{-1})$.

**Remark 3.20.** The proof of this theorem will occupy the remainder of this section. Note that the case when $\alpha = 0$, and thus $\gamma = 1$, is particularly easy, since we have

$$\chi^\lambda(h\sigma^0) = \chi^\lambda_{H_{r,n}}(h),$$

and these values are known by Theorem 2.17.

**$\kappa$-Laced Tableaux.**

Let the notations be as in Theorem 3.19, and let $\kappa = \alpha f_\lambda$. Note that the orbit of a box in $\lambda$ under the action of $\sigma^\kappa$ is of size $\gamma$.

Let $w_1$ be the permutation given in cycle notation by

$$w_1 = (1, 2, \ldots, \gamma - 1, \gamma)(\gamma + 1, \gamma + 2, \ldots, 2\gamma) \cdots . \qquad (3.21)$$

Define

$$\mathcal{L}(\lambda)^\kappa = \{L \in \mathcal{L} \mid \sigma^{-\kappa} L = w_1 L\}.$$

The elements of $\mathcal{L}(\lambda)^\kappa$ will be called $\kappa$-*laced tableaux* of shape $\lambda$. It follows from (3.5) that if $L$ is a $\kappa$-laced tableau and $1 \leq j \leq n$ then

$$\text{ct}(L(j)) = \varepsilon^{-(m\gamma - j)\kappa} \text{ct}(L(m\gamma)), \qquad (3.22)$$

where $m$ is the positive integer such that $0 \leq m\gamma - j \leq \gamma - 1$.

As an example, the necklace in Figure 3.23 is part of a standard tableau that is 3-laced. (This is the analogue of the alternating tableaux defined in [HR]).

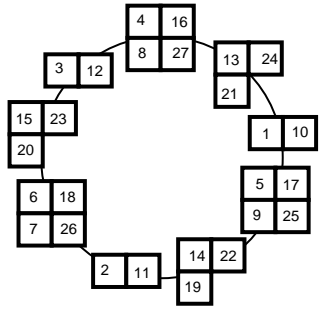

Figure 3.23. A necklace in a 3-laced standard tableau.

**Lemma 3.24.** *Let the notations be as in Theorem 3.19. If $hv_L\big|_{v_{\sigma^{-\kappa}L}} \neq 0$, then*
(a) *$L$ is $\kappa$-laced and*
(b) *every $\ell_i$ in the sequence $(\ell_1, \ldots, \ell_m)$ is divisible by $\gamma$.*

*Proof.* (a) Because of the special form of $h$ the basis elements that appear in $hv_L$ are of the form $v_{wL}$, where $w = s_{j_1} \cdots s_{j_k}$ is a product of $s_j$ such that $j_1 < j_2 < \cdots < j_k$ is a subset of the sequence $\{2, 3, \ldots, \ell_1, \ell_1 + 2, \ell_1 + 3, \ldots, \ell_2, \ell_2 + 2, \ldots\}$.

This means that, in cycle notation, $w$ is a product of cycles of the form $(i, i+1, i+2, \ldots, j-1, j)$. Thus, $hv_L\big|_{v_{\sigma^{-\kappa}L}} \neq 0$ only if $\sigma^{-\kappa}L = wL$ for some permutation of this form. But any permutation $\pi$ such that $\pi L = \sigma^{-\kappa}L$ must have all cycles of length $\gamma$, it follows that $w = w_1$ as given in (3.21). Thus, if $hv_L\big|_{v_{\sigma^{-\kappa}L}} \neq 0$, then $w_1 L = \sigma^{-\kappa}L$ and so $L$ is $\kappa$-laced.

(b) By the proof of (a), $w_1 = s_{j_1} \cdots s_{j_k}$ where $j_1 < j_2 < \cdots < j_k$ is a subset of the sequence of factors $\{2, 3, \ldots, \ell_1, \ell_1 + 2, \ell_1 + 3, \ldots, \ell_2, \ell_2 + 2, \ldots\}$. This fact and the explicit form of $w_1$ in (3.21) implies that each $\ell_i$ must be divisible by $\gamma$. $\square$

**"Dividing by $\gamma$".**

Keeping the notations as in Theorem 3.19, let us now assume that the sequence $(\ell_1, \cdots, \ell_m)$ is such that $\ell_i$ *is divisible by $\gamma$* for all $i$.

Let
$$w_j = \prod_{\substack{i \geq j \\ \gamma \nmid (i-1)}} s_i, \qquad (3.25)$$

where the product is taken with the $s_i$ in increasing order and over all $i \geq j$ such that $i - 1$ is not divisible by $\gamma$. Note that with this definition $w_1$ is the same as given in (3.21) and that $w_1 L = \sigma^{-\kappa}L$ if $hv_L\big|_{v_{\sigma^{-\kappa}L}} \neq 0$. In fact, it follows from the explicit form of $h$ and the definition of the action in (3.3) that

$$hv_L\big|_{v_{\sigma^{-\kappa}L}} = \prod_{1 \leq j \leq n} F_j(L), \qquad (3.26)$$

where $F_j(L)$ is defined as follows:
(a) $F_j(L) = (T_j)_{w_j L, w_{j+1} L}$, if $j - 1$ is not divisible by $\gamma$,
(b) $F_j(L) = (T_j)_{w_{j+1} L, w_{j+1} L}$, if $j - 1$ is divisible by $\gamma$ but $j - 1 \neq \ell_k$ for any $1 \leq k \leq m$,

(c) $F_j(L) = (t_j)^{i_k}_{w_{j+1}L, w_{j+1}L}$, if $j - 1 = \ell_k$ for some $1 \leq k \leq m$.

We shall compute the values of the $F_j(L)$ explicitly in Lemma 3.30, but first we must introduce a bit more notation.

Recall the definitions of $\bar{n}$, $\bar{r}$, $\bar{p}$, $(\bar{\ell}_1, \ldots, \bar{\ell}_m)$, $\bar{\lambda}$, and $\bar{h}$ in Theorem 3.19. If $L$ is a $\kappa$-laced tableau define integers $\rho_1, \ldots, \rho_{\bar{n}}$ and a $(d, \bar{p})$ standard tableau $\bar{L}$ as follows:

If $L(m\gamma)$ is in position $(i,j)$ of the partition $\lambda^{(k, \rho_m \bar{p} + \tau_m)}$, then
$\bar{L}(m)$ is in position $(i,j)$ of the partition $\bar{\lambda}^{(k, \tau_m)}$.

In the above $\rho_m$ and $\tau_m$ are chosen such that $0 \leq \tau_m \leq \bar{p} - 1$.

The map
$$L \longmapsto (\rho_1, \ldots, \rho_{\bar{n}}, \bar{L})$$
is a bijection between $\kappa$-laced tableaux $L$ and sequences $(\rho_1, \ldots, \rho_{\bar{n}}, \bar{L})$ where $0 \leq \rho_m \leq \gamma - 1$ for each $1 \leq m \leq n/\gamma$, and $\bar{L}$ is a $(d, \bar{p})$-standard tableau. The following is the necklace of Figure 3.23 after dividing by $\gamma = 3$.

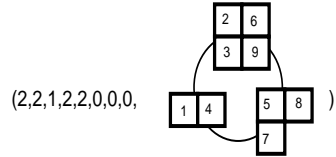

Figure 3.27. The 3-laced necklace of Figure 3.23 after division by 3.

For each $1 \leq j \leq \bar{n}$, define
$$d_j = \begin{cases} \rho_j \bmod \gamma, & \text{if } j - 1 = \bar{\ell}_k \text{ for some } k; \\ \rho_{j-1} - \rho_j \bmod \gamma, & \text{otherwise}; \end{cases}$$
(an invertible linear transformation of $(\mathbb{Z}/\gamma\mathbb{Z})^{\bar{n}}$). Then the map
$$L \longmapsto (\rho_1, \ldots, \rho_{\bar{n}}, \bar{L}) \longmapsto (d_1, \ldots, d_{\bar{n}}, \bar{L}),$$
is a bijection between $\kappa$-laced tableaux $L$ and sequences $(d_1, \ldots, d_{\bar{n}}, \bar{L})$ where $0 \leq d_m \leq \gamma - 1$ for each $1 \leq m \leq \bar{n}$, and $\bar{L}$ is a $(d, \bar{p})$-standard tableau.

The reason for introducing these bijections will become more clear in the proof of the following lemma. First let us define
$$\text{ct}(\bar{L}(m)) = \varepsilon^{\tau_m} x_k^{1/p} q^{2(j-i)}, \tag{3.28}$$
if $m$ is in position $(i,j)$ of the partition $\bar{\lambda}^{(k, \tau_m)}$ of $\bar{L}$ and then note that
$$\text{ct}(L(m\gamma)) = \varepsilon^{\rho_m \bar{p}} \text{ct}(\bar{L}(m)) = \omega^{\rho_m} \text{ct}(\bar{L}(m)), \tag{3.29}$$
where $\omega = \varepsilon^{\bar{p}}$ is a primitive $\gamma$-th root of unity.

**Lemma 3.30.** *Let the notations be as given in Theorem 3.19 and assume that $\kappa = \alpha f_\lambda$, that $L$ is a $\kappa$-laced standard tableau, and that the sequence $(\ell_1, \ldots, \ell_m)$ is such that $\ell_i$ is divisible by $\gamma$ for all $i$. Let $\omega = \varepsilon^{\bar{p}} = e^{2\pi i/\gamma}$, and let $F_j(L)$ denote the factor defined in (3.26). Let $1 \leq j \leq n$ and suppose that $k$ is such that $(k-1)\gamma < j \leq k\gamma$.*

(a) If $j-1$ is not divisible by $\gamma$, then
$$F_j(L) = (T_j)_{w_j L, w_{j+1} L} = \frac{q}{1-\varepsilon^{-(k\gamma-j+1)\kappa}} + \frac{q^{-1}}{1-\varepsilon^{(k\gamma-j+1)\kappa}}.$$

(b) If $j-1$ is divisible by $\gamma$ but $j-1 \neq \ell_i$ for any $1 \leq i \leq m$,
$$F_j(L) = (T_j)_{w_{j+1}L, w_{j+1}L} = \frac{q-q^{-1}}{1-\omega^{d_k}\frac{\mathrm{ct}(\bar{L}((k-1)))}{\mathrm{ct}(\bar{L}(k))}}.$$

(c) If $j-1 = \ell_{i-1} = \bar{\ell}_{i-1}\gamma$ for some $0 \leq i \leq m$, then
$$F_j(L) = (t_j)^{i_k}_{w_{j+1}L, w_{j+1}L} = \omega^{d_{\bar{\ell}_i-1}+1}\mathrm{ct}(\bar{L}(\bar{\ell}_{i-1}+1))^{i_k}.$$

*Proof.* (a) If $j-1$ is not divisible by $\gamma$ then $w_{j+1}L(j-1) = L(j-1)$ and $w_{j+1}L(j) = L(k\gamma)$. Thus
$$\mathrm{ct}(w_{j+1}L(j-1)) = \mathrm{ct}(L(j-1)) = \varepsilon^{-\kappa(k\gamma-(j-1))}\mathrm{ct}(L(k\gamma)), \quad \text{and}$$
$$\mathrm{ct}(w_{j+1}L(j)) = \mathrm{ct}(L(k\gamma)).$$
It follows that
$$F_j(L) = (T_j)_{w_j L, w_{j+1}L} = q^{-1} + \frac{q-q^{-1}}{1-\frac{\mathrm{ct}(w_{j+1}L(j-1))}{\mathrm{ct}(w_{j+1}L(j))}}$$
$$= q^{-1} + \frac{q-q^{-1}}{1-\varepsilon^{-\kappa(k\gamma-j+1)}}$$
$$= \frac{q}{1-\varepsilon^{-(k\gamma-j+1)\kappa}} + \frac{q^{-1}}{1-\varepsilon^{(k\gamma-j+1)\kappa}}.$$

(b) If $(k-1)\gamma = j-1$ and $j-1 \neq \ell_i$ then $w_{j+1}L(j-1) = L((k-1)\gamma)$ and $w_{j+1}L(j) = L(k\gamma)$, and
$$\mathrm{ct}(w_{j+1}L(j-1)) = \mathrm{ct}(L((k-1)\gamma)) \quad \text{and} \quad \mathrm{ct}(w_{j+1}L(j)) = \mathrm{ct}(L(k\gamma)).$$
Thus
$$F_j(L) = (T_j)_{w_{j+1}L, w_{j+1}L} = \frac{q-q^{-1}}{1-\frac{\mathrm{ct}(L((k-1)\gamma))}{\mathrm{ct}(L(k\gamma))}}$$
$$= \frac{q-q^{-1}}{1-\omega^{\rho_{k-1}-\rho_k}\frac{\mathrm{ct}(\bar{L}((k-1)))}{\mathrm{ct}(\bar{L}(k))}} = \frac{q-q^{-1}}{1-\omega^{d_k}\frac{\mathrm{ct}(\bar{L}((k-1)))}{\mathrm{ct}(\bar{L}(k))}}.$$

(c) If $j-1 = \ell_{i-1} = \bar{\ell}_{i-1}\gamma$ then $w_{j+1}L(j) = L(\ell_{i-1} + \gamma) = L((\bar{\ell}_{i-1}+1)\gamma)$ and
$$\mathrm{ct}(w_{j+1}L(j)) = \mathrm{ct}(L((\bar{\ell}_{i-1}+1)\gamma)).$$
Thus,
$$F_j(L) = (t_j)^{i_k}_{w_{j+1}L, w_{j+1}L} = \mathrm{ct}(L((\bar{\ell}_{i-1}+1)\gamma))^{i_k}$$
$$= \omega^{\rho_{\bar{\ell}_i-1}+1}\mathrm{ct}(\bar{L}(\bar{\ell}_{i-1}+1))^{i_k} = \omega^{d_{\bar{\ell}_i-1}+1}\mathrm{ct}(\bar{L}(\bar{\ell}_{i-1}+1))^{i_k}. \quad \square$$

Note that the product of the factors of type (a) in the previous lemma satisfy
$$C = \prod_{\gamma \nmid (j-1)} F_j(L) = \prod_{i=1}^{\gamma}\left(\frac{q}{1-\varepsilon^{-i}} + \frac{q^{-1}}{1-\varepsilon^i}\right)^{\bar{n}}, \tag{3.31}$$

and thus we have that

$$hv_L\big|_{v_{\sigma^{-\kappa}L}} = \prod_{1\leq j\leq n} F_j(L) = C\prod_{1\leq k\leq \bar{n}} \bar{F}_k(\bar{L}), \qquad (3.32)$$

where we define $\bar{F}_k(\bar{L}) = F_{(k-1)\gamma+1}(L)$. With this notation, the only factor in (3.32) which depends on the number $d_i$ is $\bar{F}_i(\bar{L})$.

**Proof of Theorem 3.19.**

*Proof.* Let $\kappa = \alpha f_\lambda$. Then

$$\chi^\lambda(h\sigma^\kappa) = \sum_{L\in\mathcal{L}(\lambda)} h\sigma^\kappa v_L\big|_{v_L} = \sum_{L\in\mathcal{L}(\lambda)} hv_{\sigma^\kappa L}\big|_{v_L} = \sum_{L\in\mathcal{L}(\lambda)^\kappa} hv_L\big|_{v_{\sigma^{-\kappa}L}}$$

$$= \sum_{\bar{L}\in\mathcal{L}(\bar\lambda)} \sum_{d_1,\ldots,d_{\bar n}=0}^{\gamma-1} \prod_j F_j(L)$$

$$= \sum_{\bar{L}\in\mathcal{L}(\bar\lambda)} \sum_{d_1,\ldots,d_{\bar n}=0}^{\gamma-1} C\prod_{k=1}^{\bar n} \bar{F}_k$$

$$= C\sum_{\bar{L}\in\mathcal{L}(\bar\lambda)} \prod_{k=1}^{\bar n} \left(\sum_{d_k=0}^{\gamma-1} \bar{F}_k(\bar{L})\right),$$

since the only factor in $\prod_{k=1}^{\bar n} \bar{F}_k(\bar{L})$ which depends on the number $d_i$ is $\bar{F}_i(\bar{L})$.

(a) It follows from Lemma 3.24, that if there is some $\ell_i$ that is not divisible by $\gamma$, then

$$\chi^\lambda(h\sigma^\kappa) = 0.$$

(b) Suppose that all $\ell_i$ are divisible by $\gamma$ and that $i_k \neq 0$ for some $1\leq k\leq m$. Let $j = \ell_{k-1} + 1$. Then

$$\sum_{d_{\bar\ell_{k-1}+1}=0}^{\gamma-1} \bar{F}_{\bar\ell_{k-1}+1} = \sum_{d_{\bar\ell_{k-1}+1}=0}^{\gamma-1} (t_j)^{i_k}_{w_{j+1}L, w_{j+1}L}$$

$$= \sum_{d_{\ell_{k-1}+1}=0}^{\gamma-1} \omega^{i_k d_{\bar\ell_{k-1}+1}} \mathrm{ct}(\bar{L}(\bar\ell_{k-1}+1))^{i_k} = 0,$$

and it follows that, if $i_k \neq 0$ for some $k$, then

$$\chi^\lambda(h\sigma^\kappa) = \sum_{\bar{L}\in\mathcal{L}^\kappa} hv_{\sigma^\kappa L}\big|_{v_L} = 0.$$

(c) Suppose that all $\ell_i$ are divisible by $\gamma$ and that all $i_k = 0$. Then

$$\sum_{d_{\bar\ell_{k-1}+1}=0}^{\gamma-1} \bar{F}_{\bar\ell_{k-1}+1} = \sum_{d_{\bar\ell_{k-1}+1}=0}^{\gamma-1} (t_j)^{i_k}_{w_{j+1}L, w_{j+1}L} = \sum_{d_{\bar\ell_{k-1}+1}=0}^{\gamma-1} 1 = \gamma$$

and
$$\sum_{d_k=0}^{\gamma-1} (T_j)_{w_{j+1}L, w_{j+1}L} = \sum_{d_k=0}^{\gamma-1} \frac{q - q^{-1}}{1 - \omega^{d_k} \frac{\mathrm{ct}(\bar{L}(k-1))}{\mathrm{ct}(\bar{L}(k))}}$$
$$= \frac{(q - q^{-1})\gamma}{1 - \frac{\mathrm{ct}(\bar{L}(k-1))^\gamma}{\mathrm{ct}(\bar{L}(k))^\gamma}}$$
$$= \frac{\gamma}{[\gamma]} \left( \frac{q^\gamma - q^{-\gamma}}{1 - \frac{\mathrm{ct}(\bar{L}(m-1))^\gamma}{\mathrm{ct}(\bar{L}(m))^\gamma}} \right)$$

where $j - 1 = (k - 1)\gamma$ and $[\gamma] = (q^\gamma - q^{-\gamma})/(q - q^{-1})$.

It follows that, if $i_k = 0$ for all $k$, then

$$\chi^\lambda(h\sigma^\kappa) = \sum_{L \in \mathcal{L}(\lambda)^\kappa} h v_{\sigma^\kappa L}\Big|_{v_L} = C\gamma^m \sum_{\bar{L} \in \bar{\mathcal{L}}(\bar{\lambda})} \prod_{\substack{1 \le k \le \bar{n} \\ k \ne \bar{\ell}_i}} \frac{\gamma}{[\gamma]} \left( \frac{q^\gamma - q^{-\gamma}}{1 - \frac{\mathrm{ct}(\bar{L}(k-1))^\gamma}{\mathrm{ct}(\bar{L}(k))^\gamma}} \right).$$

With the definitions of $H_{\bar{r},\bar{n}}$ as in the statement of the theorem, (2.11) and Proposition (2.12) imply that

$$\sum_{\bar{L} \in \bar{\mathcal{L}}(\bar{\lambda})} \prod_{\substack{1 \le k \le \bar{n} \\ k \ne \bar{\ell}_i}} \left( \frac{q^\gamma - q^{-\gamma}}{1 - \frac{\mathrm{ct}(\bar{L}(k-1))^\gamma}{\mathrm{ct}(\bar{L}(k))^\gamma}} \right) = \chi^{\bar{\lambda}}_{H_{\bar{r},\bar{n}}}(\bar{h}).$$

Thus,
$$\chi^\lambda(h\sigma^{\alpha f}) = C \frac{\gamma^{\bar{n}}}{[\gamma]^{\bar{n}-m}} \chi^{\bar{\lambda}}_{H_{\bar{r},\bar{n}}}(\bar{h}),$$

where $H_{\bar{r},\bar{n}}$ is as in the statement of the theorem. $\square$

## 4. THE POSET THEOREM

Curtis Greene [Gre] uses the theory of partially ordered sets (posets) and Möbius functions to prove a rational function identity ([Gre], Theorem 3.3) which can be used to derive the Murnaghan-Nakayama rule for symmetric group characters. In [HR], we modify Greene's theorem so that it can be applied to computing Murnaghan-Nakayama rules for the irreducible characters of the Iwahori-Hecke algebras of type $A_{n-1}, B_n$, and $D_n$. In this section, we extend the poset theorem of [HR] so that it can be applied to computing Murnaghan-Nakayama rules (Theorem 2.17) for the irreducible characters of the cyclotomic Iwahori-Hecke algebras of type B.

A poset is *planar* in the (strong) sense if its Hasse diagram may be order-embedded in $\mathbb{R} \times \mathbb{R}$ without edge crossings even when extra bottom and top elements are added (see [Gre] for details). A linear extension of a poset $P$ is a poset $L$ with the same underlying set as $P$ and such that the relations in $L$ form an extension of the relations in $P$ to a total order. We will denote by $\mathcal{L}(P)$ the set of all linear extensions $L$ of $P$.

The *Möbius function* of a poset $P$ is the function $\mu : P \times P \to \mathbb{Z}$ defined inductively for elements $a, b \in P$ by

$$\mu(a,b) = \mu_P(a,b) = \begin{cases} 1 & \text{if } a = b, \\ -\sum_{a \leq x < b} \mu(a,x) & \text{if } a < b, \\ 0 & \text{if } a \not\leq b. \end{cases} \quad (4.1)$$

(See [Sta] for more details on Möbius functions).

Throughout this section, $\hat{P}$ will denote a planar poset with unique minimal element $u$, and $P = \hat{P} - \{u\}$ will be the poset obtained by removing the minimal element $u$ from $\hat{P}$. We let SC be the set of minimal elements of $P$, and we call these elements *sharp corners*. Two sharp corners $s_1$ and $s_2$ of SC are "adjacent" if they are not separated by another sharp corner as the boundary of $P$ is traversed. If $s_1$ and $s_2$ are adjacent elements of SC and the least common multiple $s_1 \vee s_2$ exists, then we call $s_1 \vee s_2$ a *dull corner* of $P$. We let DC denote the set of all dull corners of $P$. Finally, we let $cc$ denote the number of connected components of $P$, and note that $cc = |\text{SC}| - |\text{DC}|$.

Let $\{x_a, a \in \hat{P}\}$, be a set of commutative variables indexed by the elements of $\hat{P}$. For each $0 \leq k \leq r-1$ and each pair $a < b$ in $\hat{P}$, define a weight, $wt^{(k)}(a,b)$, by

$$\begin{aligned} wt^{(k)}(a,b) &= \frac{1 - x_a x_b^{-1}}{q - q^{-1}} \quad \text{for all } a, b \in P, \text{ and} \\ wt^{(k)}(u,a) &= x_a^{-k} \quad \text{for all } a \in P. \end{aligned} \quad (4.2)$$

Then for any planar poset $\hat{P}$ with unique minimal element $u$, define

$$\Delta^{(k)}(\hat{P}) = \prod_{\substack{a,b \in \hat{P} \\ a \neq b}} wt^{(k)}(a,b)^{\mu_{\hat{P}}(a,b)}, \quad (4.3)$$

where $\mu_{\hat{P}}(a,b)$ is the Möbius function for the poset $\hat{P}$.

In [HR], Theorem 5.3, it is proved that

$$\sum_{\hat{L} \in \mathcal{L}(\hat{P})} \Delta^{(0)}(\hat{L}) = \Delta^{(0)}(P)(q - q^{-1})^{cc-1}, \quad (4.4)$$

and

$$\sum_{\hat{L}\in\mathcal{L}(\hat{P})} \Delta^{(1)}(\hat{L}) = \Delta^{(1)}(P) 0^{cc-1} \left(\prod_{s\in\text{SC}} x_s\right)\left(\prod_{d\in\text{DC}} x_d^{-1}\right), \quad (4.5)$$

The expansion in (4.5) is equal to zero if there is more than one connected component in $P$.

The following is our extension of the poset theorem to include values of $k > 1$.

**Theorem 4.6.** *Let $\hat{P}$ be a planar poset (as defined above) with unique minimal element $u$. Let $P = \hat{P} \setminus \{u\}$. Then*

$$\sum_{\hat{L}\in\mathcal{L}(\hat{P})} \Delta^{(0)}(\hat{L}) = (q - q^{-1})^{cc-1} \Delta^{(0)}(P),$$

*and, for $1 \leq k \leq r - 1$,*

$$\sum_{\hat{L}\in\mathcal{L}(\hat{P})} \Delta^{(k)}(\hat{L}) = \Delta^{(k)}(P)(-q + q^{-1})^{cc-1} \left(\prod_{s\in\text{SC}} x_s\right)\left(\prod_{d\in\text{DC}} x_d^{-1}\right)$$

$$\times \sum_{t=0}^{|\text{DC}|} (-1)^t e_t(x_{\text{DC}}) h_{k-t-cc}(x_{\text{SC}})$$

*where $cc$ is the number of connected components of $P$, $e_t(x_{\text{DC}})$ is the elementary symmetric function in the variables $\{x_d, d \in \text{DC}\}$, and $h_{k-t-cc}(x_{\text{SC}})$ is the homogeneous symmetric function in the variables $\{x_s, s \in \text{SC}\}$.*

*Proof.* For each $s \in \text{SC}$ define $\hat{P}_s$ to be the same poset as $\hat{P}$ except with the additional relations $s \leq s'$, for $s \neq s' \in \text{SC}$, and all other relations implied by transitivity. Each poset $P_s$ is planar, and each linear extension of $P_s$ must place the sharp corner $s$ (a minimal element of $P$) immediately after $u$ in the ordering, so we have

$$\sum_{\hat{L}\in\mathcal{L}(\hat{P})} \Delta^{(k)}(\hat{L}) = \sum_{s\in\text{SC}} \sum_{\hat{L}_s\in\mathcal{L}(\hat{P}_s)} \Delta^{(k)}(\hat{L}_s) = \sum_{s\in\text{SC}} wt^{(k)}(u,s)^{-1} \sum_{L_s\in\mathcal{L}(P_s)} \Delta^{(k)}(L_s),$$

In $P$ we have $wt^{(k)}(a,b) = wt^{(0)}(a,b)$, so by [HR], Theorem 5.3, the second sum can be computed as

$$\sum_{L_s\in\mathcal{L}(P_s)} \Delta^{(k)}(L_s) = \sum_{L_s\in\mathcal{L}(P_s)} \Delta^{(0)}(L_s) = \Delta^{(0)} = \Delta^{(k)}(P_s),$$

since $P_s$ is connected. Moreover, $wt^{(0)}(u,s)^{-1} = 1$ for each $s \in \text{SC}$, so the case $k = 0$ is proved.

From now on assume that $1 \leq k \leq r-1$. Then we have

$$\sum_{\hat{L} \in \mathcal{L}(\hat{P})} \Delta^{(k)}(\hat{L}) = \sum_{s \in \text{SC}} wt^{(k)}(u,s)^{-1} \Delta^{(k)}(P_s)$$

$$= \Delta^{(k)}(P) \sum_{s \in \text{SC}} wt^{(k)}(u,s)^{-1} \frac{\Delta^{(k)}(P_s)}{\Delta^{(k)}(P)}$$

$$= \Delta^{(k)}(P) \sum_{s \in \text{SC}} wt^{(k)}(u,s)^{-1} \prod_{\substack{a,b \in P \\ a \neq b}} \frac{wt^{(k)}(a,b)^{\mu_{P_s}(a,b)}}{wt^{(k)}(a,b)^{\mu_P(a,b)}} \quad (4.7)$$

$$= \Delta^{(k)}(P) \sum_{s \in \text{SC}} wt^{(k)}(u,s)^{-1} \prod_{\substack{a,b \in P \\ a \neq b}} wt^{(k)}(a,b)^{\mu_{P_s}(a,b) - \mu_P(a,b)}.$$

where $\mu_{P_i}(a,b)$ and $\mu_P(a,b)$ are the Möbius functions for their respective posets.

We use the work of Greene [Gre] to compute the differences $\mu_{P_s}(a,b) - \mu_P(a,b)$ for $a < b \in P$. Let $P^*$ and $P_s^*$ denote the dual of $P$ and $P_s$, respectively (that is $u \leq_{P^*} v \Leftrightarrow v \leq_P u$). Then, $\mu_P(u,v) = \mu_{P^*}(v,u)$ (see [Sta], p. 120), so we want to compute

$$\mu_{P_s^*}(b,a) - \mu_{P^*}(b,a) \quad \text{for } a < b \in P.$$

Using the Möbius notation of [Gre] (p. 8, formulas (7) and (8)), let

$$\delta_a = \sum_{t \leq_{P^*} a} \mu_{P^*}(t,a) t, \quad \text{and} \quad \delta_a^{(s)} = \sum_{t \leq_{P_s^*} a} \mu_{P_s^*}(t,a) t,$$

so that

$$a = \sum_{t \leq_{P^*} a} \delta_t, \quad \text{and} \quad a = \sum_{t \leq_{P_s^*} a} \delta_t^{(s)}.$$

In this way, $\delta_a^{(s)} = \delta_a$ for all $a \in P^* \setminus \text{SC}$, and

$$\delta_s^{(s)} = \delta_s - \sum_{s' \in \text{SC} \setminus \{s\}} s' + \sum_{d \in \text{DC}} d.$$

It follows that

$$\mu_{P_s^*}(s',s) - \mu_{P^*}(s',s) = -1, \quad \text{for all } s' \in \text{SC} \setminus \{s\},$$
$$\mu_{P_s^*}(d,s) - \mu_{P^*}(d,s) = +1, \quad \text{for all } d \in \text{DC},$$
$$\mu_{P_s^*}(a,b) - \mu_{P^*}(a,b) = 0, \quad \text{for all other } a,b \in P,$$

and

$$\prod_{\substack{a,b \in P \\ a \neq b}} wt^{(k)}(a,b)^{\mu_{P_s}(a,b) - \mu_P(a,b)} = \prod_{s' \in \text{SC} \setminus \{s\}} wt^{(k)}(s,s')^{-1} \prod_{d \in \text{DC}} wt^{(k)}(s,d).$$

Substituting back into (4.7) gives

$$\sum_{\hat{L} \in \mathcal{L}(\hat{P})} \Delta^{(k)}(\hat{L}) = \Delta^{(k)}(P) \sum_{s \in \text{SC}} wt^{(k)}(u,s)^{-1} \prod_{s' \in \text{SC} \setminus \{s\}} wt^{(k)}(s,s')^{-1} \prod_{d \in \text{DC}} wt^{(k)}(s,d).$$

Using the fact that $|SC| - |DC| = cc$ (the number of connected components of $p$), we cancel factors of $q - q^{-1}$ and factor out $x_s$ and $x_d^{-1}$ as follows

$$\sum_{\hat{L} \in \mathcal{L}(\hat{P})} \Delta^{(k)}(\hat{L}) = \Delta^{(k)}(P) \sum_{s \in SC} x_s^k \left( \prod_{s' \in SC \setminus \{s\}} \frac{q - q^{-1}}{1 - x_s x_{s'}^{-1}} \right) \left( \prod_{d \in DC} \frac{1 - x_s x_d^{-1}}{q - q^{-1}} \right)$$

$$= \Delta^{(k)}(P) \left( \prod_{s \in SC} x_s \right) \left( \prod_{d \in DC} x_d^{-1} \right) (q - q^{-1})^{cc-1} \sum_{s \in SC} x_s^{k-1} \frac{\prod_{d \in DC}(x_d - x_s)}{\prod_{s' \in SC \setminus \{s\}}(x_{s'} - x_s)}.$$

For notational convenience, let

$$F = \Delta^{(k)}(P) \left( \prod_{s \in SC} x_s \right) \left( \prod_{d \in DC} x_d^{-1} \right) (q - q^{-1})^{cc-1}.$$

Order the sharp corners $s_1, s_2, \ldots, s_{|SC|}$ from left to right as the boundary of $P$ is traversed, and let $|s_i| = i$ (its position in the ordering). Then

$$\sum_{\hat{L} \in \mathcal{L}(\hat{P})} \Delta^{(k)}(\hat{L}) = F \sum_{s \in SC} x_s^{k-1} \frac{\prod_{d \in DC}(x_d - x_s)}{\prod_{s' \in SC \setminus \{s\}}(x_{s'} - x_s)}$$

$$= F \sum_{s \in SC} x_s^{k-1} (-1)^{|SC|-|s|} \frac{\prod_{d \in DC}(x_d - x_s) \prod_{p < q \in SC \setminus \{s\}}(x_p - x_q)}{\prod_{s' < s}(x_{s'} - x_s) \prod_{s' > s}(x_s - x_{s'}) \prod_{p < q \in SC \setminus \{s\}}(x_p - x_q)}$$

$$= F \sum_{s \in SC} x_s^{k-1} (-1)^{|SC|-|s|} \frac{\prod_{d \in DC}(x_d - x_s) \prod_{p < q \in SC \setminus \{s\}}(x_p - x_q)}{V(x_{SC})},$$

where $V(x_{SC})$ is the Vandermonde determinant in the variables $\{x_s, s \in SC\}$. Moreover,

$$\prod_{d \in DC}(x_d - x_s) = \sum_{t=0}^{|DC|} e_t(x_{DC})(-1)^{|DC|-t} x_s^{|DC|-t},$$

where $e_t(x_{DC})$ is the elementary symmetric function in the variables $\{x_d, d \in DC\}$.

Again for notational convenience, let

$$G = \prod_{p < q \in SC \setminus \{s\}}(x_p - x_q).$$

We then have

$$\sum_{\hat{L} \in \mathcal{L}(\hat{P})} \Delta^{(k)}(\hat{L}) = F \frac{\displaystyle\sum_{s \in \text{SC}} x_s^{k-1}(-1)^{|\text{SC}|-|s|} \sum_{t=0}^{|\text{DC}|} e_t(x_{\text{DC}})(-1)^{|\text{DC}|-t} x_s^{|\text{DC}|-t} G}{V(x_{\text{SC}})}$$

$$= F \sum_{t=0}^{|\text{DC}|} (-1)^t e_t(x_{\text{DC}}) \frac{\displaystyle\sum_{s \in \text{SC}} x_s^{|\text{DC}|-t+k-1}(-1)^{|\text{DC}|+|\text{SC}|-|s|} G}{V(x_{\text{SC}})}$$

$$= F(-1)^{|\text{SC}|-|\text{DC}|-1} \sum_{t=0}^{|\text{DC}|} (-1)^t e_t(x_{\text{DC}}) \frac{\displaystyle\sum_{s \in \text{SC}} x_s^{|\text{DC}|-t+k-1}(-1)^{|s|-1} G}{V(x_{\text{SC}})}.$$

Notice that $x_s^{|\text{DC}|-t+k-1} = x_s^{k-t-(|\text{SC}|-|\text{DC}|)+|\text{SC}|-1} = x_s^{k-t-cc+|\text{SC}|-1}$ and that the numerator

$$\sum_{s \in \text{SC}} x_s^{k-t-cc+|\text{SC}|-1}(-1)^{|s|-1} \prod_{p<q \in \text{SC}\setminus\{s\}} (x_p - x_q) \tag{4.8}$$

is the alternating symmetrization of the monomial

$$x_{s_1}^{k-t-cc+|\text{SC}|-1} x_{s_2}^{|\text{SC}|-2} x_{s_3}^{|\text{SC}|-3} \cdots x_{s_{|\text{SC}|}}^{|\text{SC}|-|\text{SC}|}.$$

When we divide the numerator (4.8) by the Vandermonde $V(x_{\text{SC}})$, we get the Schur function $s_{(k-t-cc,0,0,\ldots,0)}(x_{\text{SC}})$ or, equivalently, the homogeneous symmetric function $h_{k-t-cc}(x_{\text{SC}})$, and so

$$\sum_{\hat{L} \in \mathcal{L}(\hat{P})} \Delta^{(k)}(\hat{L}) = F(-1)^{cc-1} \sum_{t=0}^{|\text{DC}|} (-1)^t e_t(x_{\text{DC}}) h_{k-t-cc}(x_{\text{SC}}),$$

and the proof is completed. $\square$

**Special Cases.**

The homogeneous symmetric function satisfies $h_m(x_{\text{SC}}) = 0$ unless $m > 0$, so if $k \geq 1$, then $h_{k-t-cc}(x_{\text{SC}}) = 0$, unless $cc \leq k$. In particular, when $k = 1$, the poset $P$ must be connected ($cc = 1$), and

$$\sum_{t=0}^{|\text{DC}|} (-1)^t e_t(x_{\text{DC}}) h_{k-t-cc}(x_{\text{SC}}) = e_0(x_{\text{DC}}) h_0(x_{\text{SC}}) = 1 \tag{4.9}$$

and

$$\sum_{\hat{L} \in \mathcal{L}(\hat{P})} \Delta^{(1)}(\hat{L}) = \Delta^{(1)}(P) \left( \prod_{s \in \text{SC}} x_s \right) \left( \prod_{d \in \text{DC}} x_d^{-1} \right), \tag{4.10}$$

which agrees with (4.5).

**Shapes and Standard Tableaux.**

Theorem 4.6 reduces the problem of computing $\sum_{\hat{L} \in \mathcal{L}(\hat{P})} \Delta^{(k)}(\hat{L})$ to computing $\Delta^{(k)}(P)$. In the case where $P$ is the poset of a (skew) shape, the product $\Delta^{(k)}(P)$ is readily computed and has been done so for $q = 1$ by Green [Gre] and for generic $q$ in [HR]. The result uses the natural extension of the theory of shapes and tableaux to the theory of partially ordered sets. (For a full treatment of this subject, see [Sta], whose notation we use here).

If $\lambda$ is a shape (possibly skew), then we construct a corresponding poset $P_\lambda$ whose Hasse diagram is given by placing a node in each box of $\lambda$ and then drawing edges connecting nodes in adjacent boxes. The order relation in this poset is so that the smallest nodes are in the upper left corners. For example,

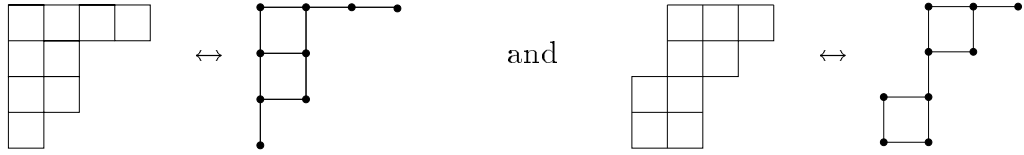

Figure 4.11

Note that posets corresponding to shapes are always planar and that the sharp and dull corners that we defined for partitions and shapes (see Figure 2.19) are exactly the sharp and dull corners of the corresponding poset.

**Theorem 4.12.** ([Gre], Theorem 3.3; [HR], Theorem 5.8) *Let $P_\lambda$ be the poset of any shape (or skew shape) $\lambda$, let $\{x_b\}$ be a set of commutative variables indexed by $\{b \in P_\lambda\}$, and let $q$ be an indeterminate. Define*

$$wt(a,b) = \frac{1 - x_a x_b^{-1}}{q - q^{-1}} \quad \text{for all } a, b \in P_\lambda \tag{4.13}$$

*and*

$$\Delta(P_\lambda) = \prod_{\substack{a,b \in P_\lambda \\ a \neq b}} wt(a,b)^{\mu_{P_\lambda}(a,b)}.$$

*Then*

$$\Delta(P_\lambda) = \left( \prod_D \frac{1 - x_b x_a^{-1}}{q - q^{-1}} \right) \left( \prod_R \frac{q - q^{-1}}{1 - x_b x_a^{-1}} \right) \left( \prod_C \frac{q - q^{-1}}{1 - x_b x_a^{-1}} \right)$$

*where*
*D is the set of pairs $(a, b)$ of boxes in $\lambda$ adjacent (northwest to southeast) in a diagonal,*
*R is the set of pairs $(a, b)$ of boxes in $\lambda$ adjacent (west to east) in a row, and*
*C is the set of pairs $(a, b)$ of boxes in $\lambda$ adjacent (north to south) in a column.*

Let $\lambda$ be a shape (or a skew shape) and let $\mathcal{L}(\lambda)$ be the set of standard tableaux of shape $\lambda$. Linear extensions of the poset $P_\lambda$ are in one-to-one correspondence with standard tableaux having skew shape $\lambda$ as follows: Given a standard tableau $T$ of shape $\lambda$ let $T(k)$ denote the box containing $k$ in $T$. Then the standard tableau $T$ corresponds to the linear extension $L$ of the poset $P_\lambda$ which has underlying set $P_\lambda$ and order relations given by $T(k) \leq_L T(l)$ if $k \leq l$. We can identify the standard tableau $T$ with the chain $L$.

Let $\hat{P}_\lambda$ be the poset $P_\lambda \cup \{u\}$ where the adjoined element $u$ satisfies $u \leq a$ for all $a \in P_\lambda$. The linear extensions of the poset $\hat{P}_\lambda$ are in one-to-one correspondence with the linear extensions of the poset $P_\lambda$. Thus, we can identify a standard tableau $T$ of shape $\lambda$ with a linear extension $\hat{L}$ of the poset $\hat{P}_\lambda$.

Let $\mu$ be the Möbius function of the linear extension $\hat{L}$ of $\hat{P}_\lambda$ that corresponds to the standard tableau $T$ of shape $\lambda$. Then, since $\hat{L}$ is a chain, $\mu$ satisfies

$$\mu(a,b) = \begin{cases} -1, & \text{if } a < b \text{ and } a \text{ is adjacent to } b \text{ in } \hat{L}, \text{ and} \\ 0, & \text{if } a < b \text{ and } a \text{ is not adjacent to } b \text{ in } \hat{L}. \end{cases}$$

It follows that

$$\Delta^{(k)}(T) = \Delta^{(k)}(\hat{L}) = \prod_{a<b\in\hat{L}} wt^{(k)}(a,b)^{\mu(a,b)} = (x_{T(1)})^k \prod_{i=2}^{n} \frac{(q-q^{-1})}{1 - x_{T(i-1)}x_{T(i)}^{-1}}.$$

**Corollary 4.14.** *Let $\lambda$ be any shape (or skew shape) with $n$ boxes. Let $\{x_b\}$ be a set of commutative variables indexed by the boxes $b \in \lambda$, and let $q$ be an indeterminate. Then*

$$\sum_{T\in\mathcal{L}(\lambda)} \Delta^{(0)}(T) = (q-q^{-1})^{cc-1} \left(\prod_D \frac{1-x_b x_a^{-1}}{q-q^{-1}}\right)\left(\prod_R \frac{q-q^{-1}}{1-x_b x_a^{-1}}\right)\left(\prod_C \frac{q-q^{-1}}{1-x_b x_a^{-1}}\right),$$

*and, for $1 \leq k \leq r-1$,*

$$\sum_{T\in\mathcal{L}(\lambda)} \Delta^{(k)}(T) = (-q+q^{-1})^{cc-1} \left(\prod_{s\in\text{SC}} x_s\right)\left(\prod_{d\in\text{DC}} x_d^{-1}\right)$$
$$\times \sum_{t=0}^{|\text{DC}|} (-1)^t e_t(x_{\text{DC}}) h_{k-t-cc}(x_{\text{SC}})$$
$$\times \left(\prod_D \frac{1-x_b x_a^{-1}}{q-q^{-1}}\right)\left(\prod_R \frac{q-q^{-1}}{1-x_b x_a^{-1}}\right)\left(\prod_C \frac{q-q^{-1}}{1-x_b x_a^{-1}}\right),$$

*where*

*cc is the number of connected components of $\lambda$,*

*SC is the set of sharp corners of $\lambda$,*

*DC is the set of dull corners of $\lambda$,*

*R is the set of pairs $(a,b)$ of boxes in $\lambda$ adjacent (west to east) in a row,*

*C is the set of pairs $(a,b)$ of boxes in $\lambda$ adjacent (north to south) in a column, and*

*D is the set of pairs $(a,b)$ of boxes in $\lambda$ adjacent (northwest to southeast) in a diagonal.*

*Proof.* Let $P_\lambda$ be the poset of the shape $\lambda$, and let $\hat{P}_\lambda$ be the poset $P_\lambda \cup \{u\}$, where the adjoined element $u$ satisfies $u \leq a$ for all $a \in P_\lambda$. Then $\hat{P}_\lambda$ is a planar poset with unique minimal element, so we apply Theorem 4.6 to compute $\sum_{T\in\mathcal{L}(\lambda)} \Delta^{(k)}(T) = \sum_{\hat{L}\in\mathcal{L}(\hat{P}_\lambda)} \Delta^{(k)}(\hat{L})$. This reduces the problem to computing $\Delta^{(k)}(P_\lambda)$. Inside of $P_\lambda$, the weights (4.2) are all independent of $k$ and of the form (4.13), so we use Theorem 4.12 to compute $\Delta^{(k)}(P_\lambda)$. □

Department of Mathematics, Macalester College, St. Paul, MN 55105
*E-mail address*: halverson@macalstr.edu

School of Mathematics and Statistics, University of Sydney, NSW 2006, Australia
*E-mail address*: ram_a@maths.su.oz.au